\documentclass[12pt]{amsart}
\usepackage{amssymb, amstext, amscd, amsmath}
\makeatletter
\def\@cite#1#2{{\m@th\upshape\bfseries%
[{#1\if@tempswa{\m@th\upshape\mdseries, #2}\fi}]}}
\makeatother
\theoremstyle{plain}
\newtheorem{thm}{Theorem}[section]
\newtheorem{cor}[thm]{Corollary}
\newtheorem{prop}[thm]{Proposition}
\newtheorem{lem}[thm]{Lemma}
\theoremstyle{definition}
\newtheorem{rem}[thm]{Remark}
\newtheorem{defn}[thm]{Definition}

\newcommand{\Prf}{\noindent\textbf{Proof.\ }}
\newcommand{\bx}{\hfill$\blacksquare$\medbreak}

\newcommand{\ca}{\mathrm{C}^*}

\newcommand{\wot}{\textsc{wot}}

\newcommand{\bbF}{{\mathbb{F}}}

\newcommand{\bbN}{{\mathbb{N}}}

\newcommand{\bbT}{{\mathbb{T}}}
\newcommand{\bbZ}{{\mathbb{Z}}}

  \newcommand{\A}{{\mathcal{A}}}
  \newcommand{\B}{{\mathcal{B}}}

  \newcommand{\E}{{\mathcal{E}}}

\renewcommand{\H}{{\mathcal{H}}}
  
  \newcommand{\J}{{\mathcal{J}}}

\renewcommand{\P}{{\mathcal{P}}}

  \newcommand{\T}{{\mathcal{T}}}

\newcommand{\fA}{{\mathfrak{A}}}

\newcommand{\fG}{{\mathfrak{G}}}

\newcommand{\qand}{\quad\text{and}\quad}

\newcommand{\qfor}{\quad\text{for}\quad}

\newcommand{\ran}{\operatorname{Ran}}

\newcommand{\spn}{\operatorname{span}}

\newcommand{\obj}{\operatorname{Obj}}

\newcommand{\fgeeplus}{\bbF^+\!(G)}
\newcommand{\fgeesplus}{\bbF^+\!(G_s)}

\begin{document}

\title[$\ca$-envelope]{The $\ca$-envelope of the tensor algebra of a directed graph}
%

\author[E. Katsoulis]{Elias~Katsoulis}
\address{Department of Mathematics\\East Carolina University\\
Greenville, NC 27858\\USA}
\email{KatsoulisE@mail.ecu.edu}
\author[D.W. Kribs]{David~W.~Kribs}
\address{Department of Mathematics and Statistics\\University of Guelph\\
Guelph, Ontario\\CANADA N1G 2W1}
\email{dkribs@uoguelph.ca}
\begin{abstract}
Given an arbitrary countable directed graph $G$ we prove the
$\ca$-envelope of the tensor algebra $\T_+(G)$ coincides with the
universal Cuntz-Krieger algebra associated with $G$. Our approach
is concrete in nature and does not rely on Hilbert module
machinery. We show how our results extend to the case of higher
rank graphs, where an analogous result is obtained for the tensor
algebra of a row-finite $k$-graph with no sources.
\end{abstract}

\thanks{Second author was partially supported by an NSERC grant}
\date{}
\maketitle

\section{Introduction}\label{S:intro}

The fundamental nonselfadjoint operator algebra associated with a
countable directed graph is its tensor algebra $\T_+(G)$
\cite{FR,MS1,KK,KP1}. Fowler, Muhly and Raeburn have recently
characterized \cite[Theorem 5.3.]{FMR} the $\ca$-envelope of the
tensor algebra of a \textit{faithful} strict Hilbert bimodule, as
the associated universal Cuntz-Pimsner algebra. When applied to
tensor algebras of graphs, their result shows that the
$\ca$-envelope of  $\T_{+}(G)$ where $G$ is a graph with
\textit{no sources}, coincides with the universal Cuntz-Krieger
algebra associated with $G$. This generalizes the well-known fact
that, for a single vertex graph with $n$ loop edges, the
$\ca$-envelope of the corresponding disc algebra is
the algebra of continuous
functions on the unit circle ($n=1$) and the Cuntz algebra ($n\geq
2$) \cite{Pop}. On the other hand, tensor algebras of graphs that do have
sources do not come from injective Hilbert bimodules and therefore
\cite[Theorem 5.3.]{FMR} does not apply to such algebras. (In that
case, Theorem 6.4 in \cite{MS2} identifies the $\ca$-envelope of
the tensor algebra of such a graph as a quotient of the associated
Toeplitz-Cuntz-Pimsner  algebra, without any further information.)
The class of tensor algebras of graphs with sources includes many
motivating examples and therefore \cite[Theorem 5.3.]{FMR} raises
the question whether or not the universal Cuntz-Krieger algebra of
an \textit{arbitrary} graph $G$ is itself the $\ca$-envelope of
$\T_{+}(G)$.

The main objective of this article is to settle this question in
the affirmative, thus providing a capstone to the investigations
from \cite[Theorem 6.8.]{Mu} and \cite{FMR, MS2, MS1} on the
$\ca$-envelope of the tensor algebra of a graph. As an added
bonus, our approach does not require any familiarity with the
heavy machinery of Hilbert bimodules. Actually, the only
`selfadjoint' prerequisite for reading this paper is the
gauge-invariant uniqueness theorem for graph $\ca$-algebras, a
rather elementary result in the theory \cite{BHRS, BPRS}.
Nevertheless, we hope the specialized techniques introduced here
will contribute to the substantial theory of tensor algebras of
Hilbert bimodules \cite{MS2,MS1}.

Another objective is to begin here a systematic study of the
$\ca$-envelope of the tensor algebra of a higher rank graph. These
were first introduced by Kumjian and Pask \cite{KumP3} as an
abstraction of the combinatorial structure underlying the higher
rank graph $\ca$-algebras of Robertson and Steger \cite{RS2,RS1}.
Their study seems in places similar to that of $1$-graphs but
there are significant complications here. One of the more notable
is the limited understanding of what it should  mean for a higher
rank graph to have a source. In this note we concentrate on the
original notion of Kumjian and Pask  of a row-finite graph with no
sources. We prove that if $\Lambda$ is such a graph, then the
$\ca$-envelope of the tensor algebra of $\T_{+}(\Lambda)$ is the
universal Cuntz-Krieger algebra $\ca(\Lambda)$.

\section{The $\ca$-envelope of $\T_{+}(G)$}

Let $G$ be a countable directed graph with vertex set $G^{0}$, edge set
$G^{1}$ and range and source maps $r$ and $s$ respectively. The \textit{Toeplitz algebra
of $G$}, denoted as $\T(G)$, is the universal $\ca$-algebra generated by a set of partial
isometries $\{ S_e\}_{ e\in G^{1}}$ and projections $\{ P_x \}_{x\in G^{0}}$ satisfying
the relations
\[
(\dagger)  \left\{
\begin{array}{lll}
(1)  & P_x P_y = 0 & \mbox{$\forall\, x,y \in G^{0}$, $ x \neq y$}  \\
(2) & S_{e}^{*}S_f = 0 & \mbox{$\forall\, e, f \in G^{1}$, $e \neq f $}  \\
(3) & S_{e}^{*}S_e = P_{s(e)} & \mbox{$\forall\, e \in G^{1}$}      \\
(4)  & \sum_{r(e)=x}\, S_e S_{e}^{*} \leq P_{x} & \mbox{$\forall\,
x \in G^{0}.$}
\end{array}
\right.
\]
The existence of such a universal object is implicit in
\cite[Theorem 3.4]{Pim} and \cite[Theorem 2.12]{MS2} and was made
explicit in \cite[Proposition 1.3 and Theorem 4.1]{FR}.

\begin{defn}\label{tensor}
Given a countable directed graph $G$, the \textit{tensor algebra of $G$},
denoted as $\T_{+}(G)$, is the norm closed subalgebra of $\T(G)$
generated by the  partial
isometries $\{ S_e\}_{ e\in G^{1}}$ and projections $\{ P_x \}_{x\in G^{0}}.$
\end{defn}

The tensor algebras associated with graphs were introduced under
the name \textit{quiver algebras} by Muhly and Solel in
\cite{Mu,MS2} as follows. Let $\lambda_{G,0}$ be the
multiplication representation of $c_0 (G^{0})$ on $l^2 (G^0)$,
determined by the counting measure on $G^0$, and let $\lambda_{G}$
denote the representation of $\T(G)$ induced by $\lambda_{G,0}$,
in the sense of \cite{MS1} and \cite{FR}. It is easily seen that
the Hilbert space  of $\lambda_G$ is $\H_G = l^{2}(\fgeeplus)$,
where $\fgeeplus$ denotes the \textit{free semigroupoid}
\cite{KP1} of the graph $G$ (also called the \textit{path space}
of $G$). This consists of all vertices $v \in G^{0}$ and all paths
$w = e_k e_{k-1} \dots e_1$, where the $e_i$ are edges satisfying
$s(e_i ) = r( e_{i-1} )$, $i = 2, 3, \dots, k$, $k \in \bbN$.
(Paths of the form $w = e_k e_{k-1} \dots e_1$ are said to have
length $k$, denoted as $|w|= k$, and vertices are called paths of
length $0$.) The maps $r$ and $s$ extend to $\fgeeplus$ in the
obvious way, two paths $w_1$ and $w_2$ are composable precisely
when $s(w_2 ) = r( w_1 )$ and, in that case, the composition $w_2
w_1$ is just the concatenation of $w_1$ and $w_2$. Let $\{ \xi_w
\}_{w \in \fgeeplus}$ denote the standard orthonormal basis of
$\H_G$, where $\xi_w$ is the characteristic function of $\{ w \}$.
Then, $\lambda_G (S_e )$, $e \in G^{1}$, is equal to the left
creation operator $L_e \in \B(\H_G)$ defined by
\[
L_e \xi_w =
\left\{
\begin{array}{ll}
\xi_{ew} & \mbox{if $ s(e) = r(w)$} \\
0 & \mbox{if $s(e)\neq r(w)$.}
\end{array}
\right.
\]
(We shall write $P_x$ for $\lambda_G (P_x)$.) Further define
$R_e$, $e\in G^1$, by the corresponding right actions on $\H_G$.
By \cite[Corollary 2.2]{FR}, the representation $\lambda_G $ is a
faithful representation of the corresponding Toeplitz algebra and
therefore a faithful representation of $\T_{+}(G)$. The algebra
$\lambda_G ( \T_{+}(G) )$ is the \textit{quiver algebra} of Muhly
and Solel \cite{MS1}. The $\wot$-closure of the quiver algebra is
the {\it free semigroupoid algebra} of the second author and Power
\cite{KP1}.

We now introduce the Cuntz-Krieger algebra of a directed graph
$G$ and one of its
faithful representations
that is convenient for the identification of the $\ca$-envelope of $\T_{+}(G)$.

Recall that a family of partial isometries $\{ S_e\}_{ e\in
G^{1}}$ and projections $\{ P_x \}_{x\in G^{0}}$ is said to obey
the Cuntz-Krieger relations associated with $G$ if and only if
they satisfy
\[
(\ddagger)  \left\{
\begin{array}{lll}
(1)  & P_x P_y = 0 & \mbox{$\forall\, x,y \in G^{0}$, $ x \neq y$}  \\
(2) & S_{e}^{*}S_f = 0 & \mbox{$\forall\, e, f \in G^{1}$, $e \neq f $}  \\
(3) & S_{e}^{*}S_e = P_{s(e)} & \mbox{$\forall\, e \in G^{1}$}      \\
(4)  & S_e S_{e}^{*} \leq P_{r(e)} & \mbox{$\forall\, e \in G^{1}$} \\
(5)  & \sum_{r(e)=x}\, S_e S_{e}^{*} = P_{x} & \mbox{$\forall\, x
\in G^{0}$ with $|r^{-1}(x)|\neq 0 , \infty$}
\end{array}
\right.
\]
The relations $(\ddagger)$ have been refined in a series of papers
by the Australian school and reached the above form in \cite{BHRS,
RS}. All refinements involved condition $(5)$ and as it stands
now, condition $(5)$ gives the equality requirement for
projections $P_{x}$ such that $x$ is not a source and receives
finitely many edges. (Indeed, otherwise condition $(5)$ would not
be a $\ca$-condition.)

It can been shown that there exists a universal $\ca$-algebra,
denoted as $\ca(G)$, associated with the relations $(\ddagger)$.
Indeed, one constructs a single family of partial isometries and
projections obeying $(\ddagger)$. Then, $\ca(G)$ is the
$\ca$-algebra generated by a `maximal' direct sum of such
families. See \cite{BPRS} for more details.

Let $G_s$ be the graph resulting by `adding tails' to $G$ in the
following sense. Let $x_1, x_2, \ldots \in G^{0}$ be the sources
of $G$. For each source $x_n$, we add to $G^{0}$ and $G^{1}$
sequences $\{ x_{n,i} \}_{i=1}^{\infty}$ and $\{ e_{n,i}
\}_{i=1}^{\infty}$ respectively, satisfying
\[
s(e_{n,i}) = x_{n,i}\quad \mbox{and}  \quad r(e_{n,i}) = x_{n,i-1}, \qquad i= 1, 2, \dots,
\]
with the convention $x_{n,0}\equiv x_n$. With these additions, the
resulting graph is denoted as $G_s$ and  it is clear that $G_s$
has no sources.

In the proof below, the faithfulness of our representation will
follow from a gauge-invariant uniqueness theorem. For that reason,
we need to construct a gauge action $\beta$ of $\bbT$ on
$\lambda_{G_s} (\T(G_s))$. For each $z \in \bbT$ we define a
unitary operator $U_{z} \in B(\H_{G_s})$ via the formula
\[U_{z}\xi_{w}= \overline{z}^{|w|}\xi_{w} \qfor w \in \fgeesplus.
\]
Note that
the family $\{ U_{z} \}_{z \in \bbT}$ induces by conjugation a gauge action  $\beta$
on the $\ca$-algebra $\lambda_{G_s} (\T(G_s))$. Indeed, if we let
\[
\beta_{z}(A) \equiv U_{z}^{*} A U_{z} \qfor A \in B(\H_{G_s}),
\]
then one easily verifies that $\beta_{z}(\lambda_{G_s}(P_x)) =
\lambda_{G_s}(P_x)$, $x \in G_s^0$, and
$\beta_{z}(\lambda_{G_s}(S_e)) = z \lambda_{G_s}(S_e)$, $e \in
G_s^1$.

\begin{thm}   \label{universal}
Let $G$ be a countable directed graph and let $G_s$ be the graph
obtained from $G$ by adding tails. Let $\ca_{\H_{G_s}}(G)
\subseteq \lambda_{G_s}(\T(G_s ))$ be the $\ca$-algebra generated
by $\{ \lambda_{G_s}(S_e) \}_{ e\in G^{1}}$ and $\{
\lambda_{G_s}(P_x) \}_{x\in G^{0}}$. Let $\pi$ be the Calkin map
on $\B(\H_{G_s})$. Then the algebra $\pi ( \ca_{\H_{G_s}}(G)) $ is
isomorphic to the universal Cuntz-Krieger $\ca$-algebra $\ca(G)$
associated with $G$.
\end{thm}

\Prf It is easy to see that the algebra $\pi ( \ca_{\H_{G_s}}(G)
)$ is generated by families of partial isometries  $\{\pi(
\lambda_{G_s}(S_e)  )\}_{ e\in G^{1}}$ and projections \break $\{
\pi(\lambda_{G_s}(P_x)) \}_{x\in G^{0}}$ satisfying $(\ddagger)$,
with respect to $G$. (By taking a quotient with the compacts, we
ensure that we get an equality in $(5)$, whenever applicable).
Since each $\lambda_{G_s}(P_x)$ has infinite dimensional range
(there are no sources in $G_s$), we obtain that $\pi
(\lambda_{G_s}(P_x) ) \neq 0$ for all $x \in G^{0}$.

Now notice that the gauge action $\beta$ is spatially implemented
and therefore passes to the Calkin algebra. It also preserves $\pi
( \ca_{\H_{G_s}}(G) )$ and therefore it induces a gauge action on
that algebra. Therefore, all the requirements of \cite[Theorem
2.1.]{BHRS} are satisfied and so $\pi ( \ca_{\H_{G_s}}(G) )$ is
isomorphic to $\ca(G)$. \bx

Theorem \ref{universal} shows that the Calkin map
\begin{equation} \label{range}
 \lambda_{G_s}(\T (G_s) ) \ni L_w \mapsto \pi(L_w) \in \pi( \lambda_{G_s}(\T( G_s )) ,
 \end{equation}
when restricted to the operators with symbols in $G$, has range
equal to the universal Cuntz-Krieger algebra for $G$. The next
lemma shows that actually (\ref{range}) is a complete isometry
when restricted to $\lambda_{G_s}(\T_{+}(G_s) )$, thus providing a
shorter proof for the main result of this section, provided that
$G$ has no sources. We note that the `unampliated' version of this
lemma has appeared as Proposition 7.3 in \cite{KK2}.

\begin{lem} \label{essnorm}
Let $\fG$ be a countable directed graph with no sources. Then,
\[ \| A \| = \|A\|_{e} , \]
for all $A\in M_n(\lambda_{\fG}(\T_{+}(\fG))$, $n \in \bbN$.
\end{lem}

\Prf It suffices to construct a sequence of isometries in the
commutant of $M_n(\lambda_{\fG}(\T_{+}(\fG))$ which converges
weakly to $0$.

Let $d\geq 1$ be a positive integer. For all $x\in\fG^0$ we may
choose a path $w_x\in \bbF^+\!(\fG)$ of length $d$ such that
$r(w_x) =x$. As the $w_x$ are distinct paths of the same length,
notice that the right creation operators $ R_{w_x}$, $x\in\fG^0$,
have pairwise orthogonal ranges. Thus, we may define an isometry
$R_{d} = \sum_{x\in \fG^{0}} R_{w_x}$ in the commutant of
$\lambda_{\fG}(\T_{+}(\fG ))$ (where the sum converges $\wot$ in
the infinite vertex case). Indeed,
\[
(R_{d})^* R_{d} = \sum_{x,y} R_{w_x}^* R_{w_y} = \sum_x
R_{w_x}^* R_{w_x} =  I.
\]

Let $R_{d}^{(n)}\equiv R_{d}\oplus R_{d} \oplus \dots \oplus R_{d}
 \in M_n(\lambda_{\fG}(\T_{+}(\fG)))^{\prime}$
denote the $n$-th inflation of $R_{d}$. Then, for an $A \in
M_n(\lambda_{\fG}(\T_{+}(\fG)))$, a compact operator $K$ and a
vector $\xi$, we have
\begin{align*}
\| A \xi \| & = \| R_{d}^{(n)} A \xi \| = \| A R_{d}^{(n)} \xi \| \\
            & \leq  \| ( A + K )R_{d}^{(n)} \xi \| + \| K R_{d}^{(n)} \xi \|  \\
            & \leq  \| A + K \| + \|K R_{d}^{(n)} \xi \|.
\end{align*}
By taking limits we obtain $  \| A \xi \| \leq  \| A + K \| $, for
all $\xi$, and so  $||A|| = ||A||_e$, as required. \bx

If $G$ has sources, in order to make use of (\ref{range}), we
identify the quiver algebra $\lambda_{G}(\T_+ (G) )$ with a
subalgebra of $\lambda_{G_s}(\T_+ (G_s) ) $.

\begin{lem}  \label{identify}
Let $G$ be a countable directed graph. Then the map
\[
\lambda_{G}(\T_+ (G) ) \ni L_w \longmapsto L_w \in
\lambda_{G_s}(\T_+ (G_s) ) , \quad w \in \fgeeplus
\]
extends to an injective $\ca$-homomorphism $ \phi : \lambda_{G}(\T
(G) ) \rightarrow \lambda_{G_s}(\T (G_s) )$.
\end{lem}

\Prf We begin by decomposing $\H_{G_s}$ into a direct sum of
reducing subspaces for $\phi(L_w)$, $w \in \fgeeplus$. First make the natural
identification of $\H_G$ inside $\H_{G_s}$. Then for all $n\geq 1$
and all $i\geq1$ (note that we do not include $i=0$ here) define
subspaces of $\H_{G_s}$ by
\[
\H_{n, i} \equiv \spn \{ \xi_{w} : w\in\bbF^+(G_s),\, s(w)=x_{n,i}
\} .
\]
Then each $\H_{n,i}$ is reducing for the range of $\phi(L_w)$, $w
\in \fgeeplus$ and
\begin{eqnarray}\label{spatialdecomp}
\H_{G_s} &=& \H_G \bigoplus \Big( \bigoplus_{n,i\geq 1}  \H_{n,i}
\Big).
\end{eqnarray}
For each $n,i$ further define $u_{n,i}\equiv e_{n,1} e_{n,2}
\cdots e_{n,i}\in\bbF^+(G_s)$ and let $\H_{n,i}^0$ be the subspace
of $\H_{n,i}$ given by
\[
\H_{n, i}^0 \equiv \spn \big\{ \xi_{wu_{n,i}} : w\in\fgeeplus,\,
s(w)=x_{n}=r(e_{n,1}) \big\} .
\]
Now let $U_{n,i}\in\B(\H_{n,i}^0,P_{x_n}\H_G)$ be the unitary
defined by
\[
U_{n,i}\, \xi_{wu_{n,i}} = \xi_w \qfor w\in\fgeeplus,\, s(w)=x_n.
\]
Then for all $n$ and $i\geq 1$ we have
\[
\phi (p(L_w ))|_{\H_{n,i}^0} = U_{n,i}^\dagger
p(L_w )U_{n,i}
\]
and
\[
\phi(p(L_w ))|_{\H_{n,i}\ominus\H_{n,i}^0} \equiv 0
\]
for any noncommutative polynomial $p(L_w)$ with $w \in \fgeeplus$.
In particular, this yields
\[
||\phi(p(L_w ))|_{\H_{n,i}}|| \leq ||
p(L_w )||.
\]
This, together with the reducing decomposition
(\ref{spatialdecomp}) of $\H_{G_s}$ and the fact that
$\phi(p(L_w ))|_{\H_G} = p(L_w )$, shows that
\begin{equation} \label{norm}
\|\phi(p(L_w )) \| =\|p(L_w ) \| ,
\end{equation}
for any noncommutative polynomial $p(L_w)$ with $w \in \fgeeplus$.
This proves the lemma since this norm estimate is obtained via a
reducing subspace decomposition. \bx

Let $\fA$ be a $\ca$-algebra and let $\A$ be a (nonselfadjoint)
subalgebra of $\fA$ which generates $\fA$ as a $\ca$-algebra and
contains a two-sided contractive approximate unit for $\fA$, i.e.,
$\A$ is an essential subalgebra for  $\fA$. A two-sided ideal $\J$
of $\fA$ is said to be a \textit{boundary ideal} for $\A$ if and
only if the quotient map $\pi: \fA \rightarrow \fA / \J$ is a
complete isometry when restricted to $\A$. It is a result of
Hamana \cite{Hamana}, following the seminal work of Arveson
\cite{Arvenv}, that there exists a boundary ideal $\J_S (\A)$, the
\textit{Shilov boundary ideal}, that contains all other boundary
ideals. In that case, the quotient $\fA / \J_S (\A)$ is called the
\textit{$\ca$-envelope} of $\A$. The $\ca$-envelope is unique in
the following sense: Assume that $\phi_1: \A \rightarrow \fA_1$ is
a completely isometric isomorphism of $\A$ onto an essential
subalgebra of a $\ca$-algebra $\fA_1$ and suppose that the Shilov
boundary for $\phi_1(\A) \subseteq \fA_1$ is zero. Then $\fA$ and
$\fA_1$ are $*$-isomorphic, via an isomorphism $\phi$ so that
$\phi( \pi(a)) = \phi_1(a)$, for all $a \in \fA$.

\begin{thm}  \label{caenvelope}
If $G$ is a countable directed graph then the $\ca$-envelope of $\T_{+}(G)$ coincides with the
universal Cuntz-Krieger algebra associated with $G$.
\end{thm}

\Prf Having (completely isometrically) identified
$\lambda_{G}(\T_{+}(G))$ with an essential subalgebra of $\ca
(G)$, the proof is now similar to that of \cite[Theorem
5.3.]{FMR}. Indeed, Theorem~\ref{universal} and Lemmata
\ref{essnorm} and \ref{identify} provide a completely isometric
isomorphism $\pi \circ \phi$ from $\lambda_{G}(\T_{+}(G))$ onto
the nonselfadjoint algebra $\ca_{+}(G)$ generated by the
generators of the universal Cuntz-Krieger algebra $\ca(G)$. In
light of the above discussion, we need to verify that the Shilov
boundary ideal $\J_S(\ca_{+}(G))$ for $\ca_{+}(G)$ inside $\ca(G)$
is zero. However, the maximality of $\J_S(\ca_{+}(G))$ and the
invariance of $\ca_{+}(G)$ under the gauge action of $\bbT$ on
$\ca(G)$ imply that $\J_S(\ca_{+}(G))$ is a gauge-invariant ideal.
Theorem 2.1 in \cite{BHRS} shows now that any non-zero
gauge-invariant ideal contains at least one of the generating
projections $P_x$, $x \in G^{0}$, or otherwise the quotient map is
an isomorphism. Hence $\J_S(\ca_{+}(G)) = \{ 0 \}$, or otherwise
the quotient map would not be faithful on $\ca_{+}(G)$. \bx

The proof above also shows the universal Cuntz-Krieger algebra of
$G$ is the $\ca$-envelope of the associated quiver algebra. To
obtain this result, we need to use the fact that $\lambda_{G}$ is
faithful \cite[Proposition 1.3 and Theorem 4.1]{FR}. Our
techniques however can give a short proof of this fact.

\begin{thm}\label{daggeruniversal}
Let $G$ be a countable directed graph.
Let $\{ P_{x}^{\prime} \}_{x \in G^0}$ and $\{
S_{e}^{\prime} \}_{e \in G^1}$ be families of projections and
partial isometries acting on a Hilbert space $\H$ which satisfy
$(\dagger)$. Then there exists a $*$-epimorphism
\[
\tau : \lambda_G(\T(G)) \longrightarrow \ca (  \{ S_{e}^{\prime}
\}_{e \in G^1} )
\]
such that $\tau(L_e) = S_{e}^{\prime}$, for all $e \in G^1$.
\end{thm}

\Prf The Wold decomposition in this setting \cite{MS1,JK1} implies
the existence of $\ca$-homomorphisms
\[
\phi_1 : \ca(G) \longrightarrow B (\H_1)
\]
and
\[
\phi_2 : \lambda_G (\T(G)) \longrightarrow B (\H_2)
\]
such that $\H = \H_1 \oplus \H_2$ and
\[
S^{\prime}_{u} = (\phi_1 \oplus \phi_2)(S_u \oplus L_u ) \qfor u
\in \fgeeplus.
\]
The desired map is then
\[
\tau = (\phi_1 \oplus \phi_2) \circ ((\pi \circ\phi) \oplus id) :
\lambda_G (\T(G)) \longrightarrow \ca ( \{ S_{e}^{\prime} \}_{e
\in \E(G)} ).
\]
Tracing the definitions of these maps shows that $\tau (L_u) =
S^{\prime}_{u}$ for all $ u \in \fgeeplus$. \bx

In \cite{FR}, Fowler and Raeburn  obtained a generalization of
Coburn's Theorem to arbitrary Toeplitz-Cuntz-Krieger
$\ca$-algebras, once again making extensive use of the Hilbert
bimodule machinery. For Toeplitz-Cuntz-Krieger $\ca$-algebras,
associated with graphs having no sources, the same result was
obtained independently in \cite{KK2} by more elementary means.
Lemma \ref{identify} allows us to apply the technique from
\cite{KK2} to the general case.

\begin{thm}   \label{Coburn}
Let $G$ be a countable directed graph. Let $\{P_{x}^{\prime}\}_{x \in G^0}$ and
$\{S_{e}^{\prime}\}_{e \in G^1}$ be families of projections and
partial isometries respectively, acting on a Hilbert space $\H$
and satisfying $(\dagger)$. If,
\begin{equation} \label{properisom}
\sum_{r(e) = x} S_{e}^{\prime} (S_{e}^{\prime})^{ *} \neq
P_{x}^{\prime}     \qfor x \in G^0,
\end{equation}
then there exists an injective $*$-homorphism $\phi: \lambda_G
(\T(G)) \rightarrow B(\H)$ such that $\phi( L_e ) =
S_{e}^{\prime}$ for all $e \in G^1$.

In particular, the $\ca$-algebra $\ca(\{S_{e}^{\prime}\}_{e \in
G^1})$, generated by the collection $\{S_{e}^{\prime}\}_{e \in
G^1}$, is isomorphic to $\lambda_G (\T(G))$.
\end{thm}

For a proof, repeat the arguments in \cite{KK2}.

\section{Higher rank graphs}\label{S:higherrank}

In this section we define the natural generalization of the tensor
algebra for a higher rank graph and, using a concrete approach,
prove its $\ca$-envelope is the corresponding universal
Cuntz-Krieger algebra when the graph is row-finite with no
sources. This also yields another concrete proof for the case of
1-graphs with no sources.

\begin{defn}\label{hrdefn}
(Kumjian and Pask \cite{KumP3}) A $k$-graph $(\Lambda, d)$
consists of a countable small category $\Lambda$, with range and
source maps $r$ and $s$ respectively, together with a functor $d:
\Lambda \rightarrow \bbZ_+^k$ satisfying the factorization
property: for every $\lambda\in\Lambda$ and $m,n\in\bbZ_+^k$ with
$d(\lambda) = m+n$, there are unique elements $\mu,\nu\in\Lambda$
such that $\lambda = \mu\nu$ and $d(\mu) = m$ and $d(\nu) =n$.
\end{defn}

By the factorization property we may identify the objects
$\obj(\Lambda)$ of $\Lambda$ with the subset $\Lambda^0 \equiv
d^{-1} (0,\ldots ,0)$. We also write $\Lambda^n$ for $d^{-1}(n)$,
$n\in \bbZ_+^k$, and put $\Lambda^n(v)=\{\lambda\in\Lambda^n
:r(\lambda)=v\}$ for $v\in\Lambda^0$. Following \cite{KumP3}, if
$\Lambda^n(v)$ is nonempty for all choices of $v\in\Lambda^0$ and
$n\in\bbN^k$, we shall say $\Lambda$ has {\it no sources}.
Further, $\Lambda$ is said to be {\it row-finite} if for all
$m\in\bbN^k$ and $v\in\Lambda^0$ the set $\Lambda^m(v)$ is finite.
Let us define a grading function $\delta : \Lambda \rightarrow
\bbN$ for $\Lambda$ by $\delta(\lambda) = n_1 + \ldots + n_k$ when
$d(\lambda) =(n_1,\ldots ,n_k)$.

It is convenient to view a higher rank $k$-graph as a directed
graph with directed edges coloured one of $k$ possible colours and
vertices identified with objects. The elements of $\Lambda^1$ are
the edges of $\Lambda$ and the factorization property describes
how different `paths' are related. See \cite{FS,KP3,KumP3,RaSiY}
for further introductory discussions on the subject.

Let $\Lambda$ be a row-finite higher rank graph with no sources. A
family $\{S_\lambda:\lambda\in\Lambda\}$ is said to satisfy the
Toeplitz-Cuntz-Krieger relations for $\Lambda$ if the following
holds:
\[
(\dagger)  \left\{
\begin{array}{lll}
(1)  & S_v S_w = 0 & \mbox{$\forall\, v,w \in \Lambda^0$, $ v \neq w$}  \\
(2) & S_{\lambda\mu} = S_\lambda S_\mu & \mbox{$\forall\,\lambda,\mu\in\Lambda$}  \\
(3) & S_{\lambda}^{*}S_\lambda = S_{s(\lambda)} & \mbox{$\forall\, \lambda\in\Lambda$}      \\
(4)  &  \sum_{\lambda\in\Lambda^n(v)}\, S_\lambda S_{\lambda}^{*}
\leq S_v & \mbox{$\forall\, v\in\Lambda^0$ and $n\in\bbN^k$.}
\end{array}
\right.
\]
If equality is satisfied in condition $(4)$ for all $v$, $n$, then
$\{S_\lambda:\lambda\in\Lambda\}$ is said to satisfy $(\ddagger)$,
the analogue of the Cuntz-Krieger relations here. The algebra
$\ca(\Lambda)$ \cite{KumP3} is the universal $\ca$-algebra
generated by a family of partial isometries satisfying
$(\ddagger)$.

The Toeplitz algebra $\T(\Lambda)$ \cite{RaSi} is the universal
$\ca$-algebra generated by a family of partial isometries
satisfying $(\dagger)$. Let us define the following nonselfadjoint
algebra as a generalization of the 1-graph case.

\begin{defn}
The norm closed  subalgebra of $\T (\Lambda)$ generated by the
family $\{ S_\lambda : \lambda\in\Lambda \} $ is denoted as
$\T_+(\Lambda)$.
\end{defn}

We will work with the following faithful representation of
$\T_+(\Lambda)$. Given a higher rank graph $\Lambda$ (row-finite
with no sources), let $\H_\Lambda$ be the Fock space Hilbert space
with orthonormal basis $\{\xi_\lambda: \lambda\in\Lambda \}$.
Define partial creation operators $L=
(L_\lambda)_{\lambda\in\Lambda}$ on $\H_\Lambda$ as follows:
\[
L_\lambda \xi_{\mu} = \left\{
\begin{array}{ll}
\xi_{\lambda\mu} & \mbox{if $ s(\lambda) = r(\mu)$} \\
0 & \mbox{if $s(\lambda)\neq r(\mu)$.}
\end{array}
\right.
\]
Further define partial isometries $\{R_\lambda : \lambda \in
\Lambda\}$ by the corresponding right actions on $\H_\Lambda$. The
$\wot$-closure of this representation is a `higher rank
semigroupoid algebra' of the second author and Power \cite{KP3}.
Faithfulness of this `Fock representation' on $\T(\Lambda)$, and
hence on $\T_+(\Lambda)$, was established by Raeburn and Sims
\cite{RaSi}.

In the following discussion let $\Lambda = (\Lambda, d)$ be a
row-finite higher rank $k$-graph with no sources. The no source
assumption on $\Lambda$ allows the following choice of infinite
paths. To each vertex $v \in \Lambda^{0}$ we associate a unique
infinite path
\[
\mu_v \equiv e_{v , 1} e_{v , 2} \cdots \, , \quad e_{v , i} \in
\Lambda^1,\quad i = 1, 2, \dots
\]
with range $v$, that is $r(e_{v , 1}) = v$. Each path $\mu_v$, $v
\in \Lambda^{0}$, is chosen so that for every $m \in \bbZ_+^k$,
there exists an $i \in \bbN$ such that $d(\mu_{v, i}) \geq m$ in
$\bbN^k$, where
\[
\mu_{v, i} \equiv e_{v , 1}e_{v , 2} \cdots e_{v , i}.
\]
For convenience put $\mu_{v,0} = v$. We will denote the collection
of all such paths as
\[
\Lambda_{red} = \big\{ \mu_{v,i} : v\in\Lambda^0, \,\, i =
0,1,2,\ldots \big\}.
\]

We shall construct a faithful representation of $\ca(\Lambda)$
that will allow us to identify the $\ca$-envelope of
$\T_+(\Lambda)$. First note that if $\lambda\in\Lambda$ then we
may define $\lambda^{-1}$ as the natural equivalence class of
products of elements $e^{-1}$, where $e\in\Lambda^1$, determined
by $\lambda$. Further, for all $\lambda,\mu\in\Lambda$ with
$s(\lambda)=s(\mu)$ we define $\lambda\mu^{-1}$ to be the
equivalence class of reduced products given by applications of the
factorization property followed by cancellations of products of
the form $ee^{-1} = r(e)$. Let
\[
\Gamma =
\Lambda\Lambda^{-1}_{red}=\{\lambda\mu^{-1}:\lambda \in\Lambda,\,
 \mu \in \Lambda_{red} \, ,
s(\lambda)=s(\mu)\},
\]
where it is understood that $\lambda\mu^{-1}$ corresponds to the
equivalence class of reduced products determined by $\lambda,\mu$
as above.

Let $\H = \ell^2(\Gamma)$ be the Hilbert space with orthonormal
basis $\{\xi_{\gamma}: \gamma\in\Gamma\}$. Consider operators $S'
= (S^\prime_\lambda)_{\lambda\in\Lambda}$ defined on $\H$ by
\[
S_\lambda^\prime \xi_{\lambda_1 \mu^{-1}} = \left\{
\begin{array}{ll}
\xi_{\lambda\lambda_1\mu^{-1}} & \mbox{if $ s(\lambda) = r(\lambda_1)$} \\
0 & \mbox{if $s(\lambda)\neq r(\lambda_1)$.}
\end{array}
\right.
\]
Let  $\tau_\Lambda$ be the map that identifies the generators of
$\ca(\Lambda)$ with the operators $S' =
(S^\prime_\lambda)_{\lambda\in\Lambda}$, so that $\tau_\Lambda
(S_\lambda) = S_\lambda^\prime$ for $\lambda\in\Lambda$.

\begin{prop}\label{repn}
The map $\tau_\Lambda: \ca(\Lambda)\rightarrow\B(\H)$ defines a
faithful representation of $\ca(\Lambda)$.
\end{prop}

\Prf It is evident that  $\{S_\lambda^\prime: \lambda\in\Lambda\}$
is a family of partial isometries that satisfy conditions $(1)$,
$(2)$ and $(3)$ of the $(\dagger)$ relations above. To see
condition $(4)$, let $v\in\Lambda^0$ and let $m\in\bbN^k$. The
range of $S_v^\prime$ is given by $\ran S_v^\prime = \spn\{
\xi_\gamma : \gamma = \lambda\mu^{-1}\in\Gamma, \,
r(\lambda)=v\}$. Consider an arbitrary element $\gamma =
\lambda\mu^{-1}_{v, i}\in\Gamma$. Choose $j > i $ so that
\[
d(e_{v , i+1}e_{v , i+2} \cdots e_{v , j} ) \geq m.
\]
Then $\lambda e_{v , i+1}e_{v , i+2} \cdots e_{v , j}$ belongs to
$\Lambda$ with $d \big( \lambda e_{v , i+1}e_{v , i+2} \cdots e_{v
, j}\big) \geq m$ as well. By the factorization property there are
$\lambda_1,\lambda_2 \in \Lambda$ such that
\[
\lambda  e_{v , i+1}e_{v , i+2} \cdots e_{v , j}
=\lambda_1\lambda_2  \qand    d(\lambda_1)=m.
\]
Since $s(\lambda_2 ) = s( e_{v , j}) $, we have that
$\xi_{\lambda_2 \mu_{v, j}^{-1}} \in \H$. Furthermore,
\begin{align*}
S_{\lambda_{1}}(\xi_{\lambda_2 \mu_{v, j}^{-1}} ) &= \xi_{(\lambda
e_{v , i+1}e_{v , i+2} \cdots e_{v , j}) \mu_{v, j}^{-1}} \\
               &= \xi_{\lambda \mu_{v, i}^{-1} } = \xi_{\gamma}
\end{align*}
and so $\xi_\gamma$ belongs to $\ran S_{\lambda_1}$. Finally, the projections
$S_\lambda S_\lambda^*$, $\lambda\in\Lambda^m$, are easily seen to
have mutually orthogonal ranges by the factorization property.
Thus we have verified condition $(4)$ and $\tau_\Lambda$ defines a
representation of $\ca(\Lambda)$ on $\H$.

To show that $\tau_\Lambda$ is faithful, we use the
gauge-invariant uniqueness theorem of \cite{KumP3}. The
universality of $\ca ( \Lambda)$ implies the existence of a
canonical action of the $k$-torus $\bbT^{k}$, called the
\textit{gauge action}, $\alpha : \bbT^k \rightarrow
{\operatorname{Aut}} \, \ca(\Lambda)$ defined for $t = (t_1, t_2,
\dots ,t_k) \in \bbT^k$ and $S_\lambda \in \ca(\Lambda)$ by
\[
\alpha_{t}(S_\lambda)= t^{d(\lambda)}S_\lambda
\]
where $t^{m}= t_{1}^{m_1} t_{2}^{m_2} \cdots t_{k}^{m_k}$ for $m =
(m_1, m_2, \dots ,m_k) \in \bbT^k$. Let $\tau : \ca(\Lambda)
\rightarrow \B$ be a representation of $\ca(\Lambda)$ so that
$\tau(S_v ) \neq 0$, for all $v \in \Lambda^{0}$. The
gauge-invariant uniqueness theorem asserts that $\tau$ is faithful
provided that there exists an action $\beta : \bbT^k \rightarrow
{\operatorname{Aut}}\, \B$ so that $\tau \circ \alpha_t = \beta_t
\circ \tau$, for all $ t \in \bbT^k$.

For each $t \in \bbT^{k}$, define a unitary operator $U_t \in \B(\H)$
as
\[
U_t\, \xi_{\lambda \mu^{-1}} \equiv t^{(d(\mu)-d(\lambda))}
\xi_{\lambda \mu^{-1}}, \quad \lambda \in \Lambda, \, \mu \in
\Lambda_{\operatorname{red}},
\]
and let $\beta_t (A) \equiv U_{t}^{*} A U_{t}$, $ A \in B(\H)$. It is easy to check that
$\tau_\Lambda \circ \alpha_t = \beta_t \circ \tau_\Lambda$, for all $ t \in \bbT^k$,
and so $\tau_\Lambda$ is faithful.
\bx

\begin{lem}\label{hrlemma}
Let $\Lambda = (\Lambda,d)$ be a row-finite higher rank graph with
no sources and let $\tau_\Lambda: \ca(\Lambda) \rightarrow \B(\H)$
be as above. If $p$ is a polynomial in the polynomial ring
$\P_\Lambda^+$ generated by $\Lambda$ then $||p(L)|| =
||\tau_\Lambda(p(S))||$.
\end{lem}

\Prf For $v\in \Lambda^0$ and $\mu\in\Lambda_{red}$ with
$r(\mu)=v$, let $\H_{v,\mu}$ be the subspace of $\H$ given by
\[
\H_{v,\mu} = \spn \{ \xi_{\lambda\mu^{-1}} : \lambda\in \Lambda,\,
s(\lambda) = s(\mu)\}.
\]
Further define for $k\geq 0$ subspaces
\[
\H_{k} = \bigoplus_{v \in \Lambda^{0}} \H_{v,\mu_{v,k}}.
\]

Then $\H_{v,\mu}$ may be naturally identified with the subspace
$\H_{s(\mu )} \equiv R_{s(\mu)} \H_\Lambda = \spn\{\xi_{\lambda} :
\lambda\in\Lambda,\, s(\lambda)=s(\mu)\}$ of the Fock space
$\H_\Lambda$ via the unitary $U_{v,\mu }:\H_{v,\mu} \rightarrow
\H_{s(\mu)}$ which is defined on standard basis vectors
$\xi_{\lambda \mu^{-1}}$ of $\H_{v,\mu}$ by $U_{v,\mu
}\,\xi_{\lambda\mu^{-1}} = \xi_{\lambda}$. The construction of
$\Gamma = \Lambda\Lambda^{-1}_{red}$ shows that $U_{v,\mu}$ is
well defined and determines a unitary operator. Note that,
\[
 \tau_{\Lambda}(S_{\lambda})\mid_{\H_{v,\mu}} =
 U_{v,\mu }^{*} \big( L_{\lambda}\mid_{R_{s(\mu)}\H_{\Lambda}}\big) U_{v,\mu }
\]
for any $\lambda \in \Lambda$.

Under these identifications, the restrictions of the $\tau_\Lambda
(S_\lambda)$ to the invariant subspace $\H_0 =
\oplus_{v\in\Lambda^0} \H_{v,\mu_{v,0}} = \oplus_{v\in\Lambda^0}
\H_{v,v}$ are jointly unitarily equivalent to the creation
operators $L_\lambda$;
\begin{equation*} \label{restriction}
\tau_\Lambda (S_\lambda)|_{\H_0} \cong L_\lambda \qfor \lambda \in
\Lambda.
\end{equation*}
Thus, if $p$ is a polynomial in the polynomial ring $\P^+_\Lambda$
generated by $\Lambda$, then
\begin{equation*} \label{one direction}
\| p(L) \| = \| p(\tau_\Lambda (S))|_{\H_0} \| \leq \| p(\tau_\Lambda (S)) \|
\end{equation*}

On the other hand, the restrictions of the $\tau_\Lambda
(S_\lambda)$ to the invariant subspaces $\H_k =
\oplus_{v\in\Lambda^0}  \H_{v,\mu_{v,k}}$, for $k\geq 0$ are
jointly unitarily equivalent to a direct sum of restrictions of
$L_\lambda$ to invariant subspaces of the form $\H_{w}$,
$w\in\Lambda^0$. Therefore,
\[
 \| p(\tau_\Lambda (S))|_{\H_k} \| \leq \| p(L) \| \qfor k\geq 0.
\]
However, since $\H$ is given by the closed span $\H= \vee_{k \geq
0} \H_k$ of this increasing chain of subspaces, the sequence of
compression operators $p(\tau_\Lambda (S))P_{\H_k}=
P_{\H_k}p(\tau_\Lambda (S))P_{\H_k}$,  $ k\geq 0$,  converges
strongly to $p(\tau_\Lambda (S))$. Hence,
\[
\| p(\tau_\Lambda (S)) \|\leq \liminf_{k \geq 0} \| p(\tau_\Lambda
(S))P_{\H_k} \| \leq \| p(L) \|,
\]
and the conclusion follows. \bx

The ampliated version of the above arguments shows that the
isometric map
\begin{eqnarray*}
\T_{+}(\Lambda) &\rightarrow& \tau_{\Lambda}(\ca(\Lambda)) \\
 p(L) &\mapsto& \tau_\Lambda(p(S))
\end{eqnarray*}
of Lemma~\ref{hrlemma} is a complete isometry and so we obtain the
following.

\begin{cor}\label{cihr}
There exists a complete isometry  $\phi: \T_{+}(\Lambda) \rightarrow \ca(\Lambda)$,
which maps generators to generators.
\end{cor}

\Prf Recall that the map $\tau_{\Lambda}$ is a faithful representation of $\ca(\Lambda)$.
The conclusion now follows from the above.
\bx

\begin{thm} \label{Lidentify}
If $\Lambda$ is a row-finite higher rank graph  with no sources,
the $\ca$-envelope of $\T_{+}(\Lambda)$ coincides with the
universal Cuntz-Krieger algebra $\ca(\Lambda)$ associated with
$\Lambda$.
\end{thm}

\Prf The proof is now identical to that of
Theorem~\ref{caenvelope}. Indeed, Corollary~\ref{cihr} provides a
completely isometric isomorphism from $\T_{+}(\Lambda)$ onto the
nonselfadjoint algebra $\ca_{+}(\Lambda)$ generated by the
generators of the universal Cuntz-Krieger algebra $\ca(\Lambda)$.
In light of the discussion above Theorem~\ref{caenvelope}, we need
to verify, once again, that the Shilov boundary ideal
$\J_S(\ca_{+}(\Lambda))$ for $\ca_{+}(\Lambda)$ inside
$\ca(\Lambda)$ is zero. However, the maximality of
$\J_S(\ca_{+}(\Lambda))$ and the invariance of $\ca_{+}(G)$ under
the gauge action of $\bbT$ on $\ca(\Lambda)$ imply that
$\J_S(\ca_{+}(\Lambda))$ is a gauge-invariant ideal. Theorem 5.2
in \cite{RaSiY} shows now that any non-zero gauge-invariant ideal
contains at least one of the generating projections $P_v$, $v \in
\Lambda^{0}$. Hence $\J_S(\ca_{+}(\Lambda)) = \{ 0 \}$, or
otherwise the quotient map would not be faithful on
$\ca_{+}(\Lambda)$. \bx

\begin{rem}
In light of Theorems \ref{daggeruniversal} and \ref{Coburn}, it is
natural to expect that similar results, with similar proofs, are
valid in the higher rank context. Indeed, such results have been
obtained by Raeburn, et al,  making use of the product Hilbert
bimodule theory. However, we would expect that a Wold
decomposition type approach could yield simpler proofs.
\end{rem}


\vspace{0.1in}

{\noindent}{\it Acknowledgement.} The research in this article was
motivated by the referee's report for \cite{KK}. We are thankful
to the referee of \cite{KK} for bringing to our attention the
results in \cite{FMR} and also pointing out to us that the case of
tensor algebras of graphs with sources was not covered by the
results in \cite{FMR}.


\end{document}